\documentclass{article}
\usepackage[cp1251]{inputenc}
\usepackage[dvips]{graphicx}
\usepackage{indentfirst}
\usepackage{amsfonts,amssymb}
\font\cyr=wncyr10

\begin{document}
\title{Grid diagrams of Lorenz links}
\author{Roman Razumovsky
\footnote{This work is supported in part by RFBR (grant no.
10-01-91056-{\cyr NCNI\_a})} } \maketitle

\begin{abstract}
In paper \cite{birman} Joan Birman and Ilya Kofman prove the
coincidence of the class of Lorenz links and the class of twisted
links. The proof in that work is algebraic. We will identify this
class in terms of grid diagrams and provide a transparent
geometric argument for Birman-Kofman's result.
\end{abstract}

Lorenz links are periodic orbits of the Lorenz "strange
attractor", which arises in physics. See \cite{old} for details.
Twisted links are a generalization of torus links: it is allowed
to twist not the whole set of strings but also some of its
subsets. We will introduce a class of oriented links by putting
restrictions on their grid diagrams and then we will show that it
coincides with the class of Lorenz links and with that of twisted
links.

\textbf{Definition 1.} A permutation $\sigma \in S_n$ is called a
shuffle if there exists k such that:
$$\sigma(1)<\sigma(2)<...<\sigma(k),$$
$$\sigma(k+1)<\sigma(k+2)<...<\sigma(n).$$

Let $\sigma$ be a shuffle such that $\sigma(i)\ne i$ for all
$i=1,...,n$. Then the closure of the permutation braid associated
to the shuffle $\sigma$ is called a Lorenz link. Lorenz links
inherit their orientation from the corresponding braids. We will
associate the following vector to $\sigma$
$$\langle\sigma(1)-1,\sigma(2)-2,...,\sigma(k)-k\rangle$$
and call it the Lorenz vector of $\sigma$.

We will denote the Lorenz link corresponding to $\sigma$ by the
associated Lorenz vector. We will also use the following short
notation for this vector:
$$\langle\underbrace{p_1,p_1,...,p_1}_{q_1},\underbrace{p_2,p_2,...,p_2}_{q_2},...,\underbrace{p_s,p_s,...,p_s}_{q_s}\rangle=\langle{p_1}^{q_1},{p_2}^{q_2},...,{p_s}^{q_s}\rangle.$$

\textbf{Example 1.} In Fig.1 the Lorenz link with Lorenz vector
$\langle3^4,5^3\rangle$ is shown.

\begin{figure}[h]
\centering {
    \includegraphics[height=4cm,width=6cm]{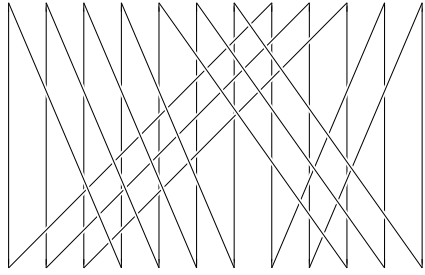}
} \caption{Lorenz link
$\langle3,3,3,3,5,5,5\rangle=\langle3^4,5^3\rangle$}
\end{figure}

\textbf{Definition 2.} Twisted links (T-links) are the closures of
braids on n strings of the following form:
$$(\sigma_{1}\sigma_{2}...\sigma_{p_1-1})^{q_1}(\sigma_{1}\sigma_{2}...\sigma_{p_2-1})^{q_2}...(\sigma_{1}\sigma_{2}...\sigma_{p_s-1})^{q_s},$$
$$p_1<p_2<...<p_s=n.$$

The T-link, associated to this braid is denoted by
T(($p_1,q_1$),...,($p_s,q_s$)) (see Fig.2). The class of T-links
is a generalization of the class of torus links, which corresponds
to $s=1$. T-links also inherit their orientations from the braids.

\textbf{Example 2.} In Fig.2 the T-link T((3,4),(5,3)) is shown.

\begin{figure}[h]
\centering {
    \includegraphics[height=4cm,width=2cm]{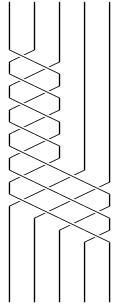}
} \caption{Braid T((3,4),(5,3))}
\end{figure}

\textbf{Definition 3.} Let $\sigma \in S_n$ be a shuffle and $v$
is its Lorenz vector. We will denote by diag($v$) (or by
diag($\sigma$)) the following link given by the coordinates of the
vertices of its grid presentation (see \cite{dynnikov} for the
definition of grid (rectangular) presentations):
$$(1,1),(2,2),...,(n,n),(1,\sigma(1)),(2,\sigma(2)),...,(n,\sigma(n)).$$

If $\sigma(i)>i$ a vertical edge with x-coordinate i is oriented
upwards; otherwise it is oriented downwards. Horizontal edges are
oriented appropriately. We will call such links diagonal.

\textbf{Example 3.} In Fig.3 the diagonal link
diag($\langle3,3,3,3,5,5,5\rangle$) is shown.

\begin{figure}[h]
\centering {
    \includegraphics[height=4cm,width=4cm]{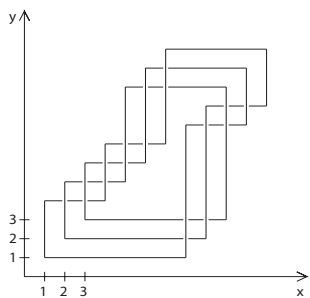}
} \caption{Link diag($\langle3,3,3,3,5,5,5\rangle$)}
\end{figure}

\textbf{Theorem 1.} The class of diagonal links coincides with the
class of Lorenz links. More precisely, the diagonal link
diag($\langle{p_1}^{q_1},{p_2}^{q_2},...,{p_s}^{q_s}\rangle$) is
isotopic to the Lorenz link
$\langle{p_1}^{q_1},{p_2}^{q_2},...,{p_s}^{q_s}\rangle$.

\textbf{Proof.} Let's denote the diagonal $x=y$ of
diag($\langle{p_1}^{q_1},{p_2}^{q_2},...,{p_s}^{q_s}\rangle$) by
$l$. Let's denote the set of vertical edges with x-coordinate less
or equal to $q_1+q_2+...+q_s$ by $L$. The remaining vertical edges
will form the set $R$.

Now we introduce the isotopy. First we rotate our grid diagram by
$45^0$ clockwise. Next we rotate each edge from $L$ by $45^0$
counterclockwise around its lower end and stretch them as
necessary so as to have all their upper ends on a straight line
parallel to $l$. Then we will apply the same procedure to $R$,
rotating the edges around their upper ends. The other edges, which
were horizontal before the rotation are adjusted appropriately.
Refer to Fig.4 for an example.

\begin{figure}[h]
\centering {
    \includegraphics[height=4cm,width=4cm]{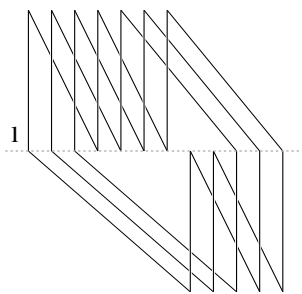}
} \caption{Link T((3,4),(5,3))}
\end{figure}

Next we fold our link along the horizontal line $l$ and put the
lower part over the upper one. Now there are only two rows of
vertices of our piecewise-linear planar diagram. There are three
layers of edges: the middle layer consists of
$q_1+q_2+...+q_s+p_s$ vertical edges, the top layer consists of
$p_s$ edges crossing the edges from the middle layer negatively,
the bottom layer consists of $q_1+q_2+...+q_s$ edges crossing the
edges from the middle layer negatively. Consider the edges of the
bottom layer. They don't cross each other, so there will be no
obstruction to pull them aside and place in the top layer. Refer
to Fig.1 for an example. This is a closed permutation braid
associated to $\sigma$, which is a shuffle. Thus it is the Lorenz
link $\langle{p_1}^{q_1},...,{p_s}^{q_s}\rangle$.

\textbf{Theorem 2.} The class of diagonal links coincides to the
class of twisted links. Precisely, the diagonal link
diag($\langle{p_1}^{q_1},{p_2}^{q_2},...,{p_s}^{q_s}\rangle$) is
isotopic to the twisted link $T((p_1,q_1),...,(p_s,q_s))$.

\textbf{Proof.} Let's cut open the horizontal edges of
diag($\langle{p_1}^{q_1},{p_2}^{q_2},...,{p_s}^{q_s}\rangle$) with
y-coordinate $\sigma(i)$ $(i>q_1+q_2+...+q_s)$ and connect them
"through the infinity". Refer to Fig.5 for an example.

\begin{figure}[h]
\centering {
    \includegraphics[height=4cm,width=4cm,angle=90]{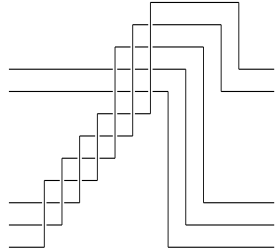}
} \caption{Link T((3,4),(5,3))}
\end{figure}

This is a diagram of a braid on $p_s$ strings, oriented
rightwards. The strings are numbered bottom-to-up. All crossings
are placed over the former diagonal of the diagonal link. They are
induced by the vertical edges of the diagonal link with
x-coordinate $i \le q_1+q_2+...+q_s$. Each such edge contributes
$\sigma_{1}\sigma_{2}...\sigma_{\sigma(i)-1}$ to the braid word
since the left string crosses over the next $\sigma(i)-1$ strings.
Finally we get the following braid word:
$$(\sigma_{1}\sigma_{2}...\sigma_{p_1-1})^{q_1}(\sigma_{1}\sigma_{2}...\sigma_{p_2-1})^{q_2}...(\sigma_{1}\sigma_{2}...\sigma_{p_s-1})^{q_s},$$

Thus we got T(($p_1,q_1$),...,($p_s,q_s$)).

\textbf{Corollary(\cite{birman}).}The class of Lorenz links
coincides to the class of twisted links. Precisely, the Lorenz
link $\langle{p_1}^{q_1},{p_2}^{q_2},...,{p_s}^{q_s}\rangle$ is
isotopic to the twisted link T(($p_1,q_1$),...,($p_s,q_s$)).

\end{document}